\documentclass{svjour3}
\usepackage[latin1]{inputenc}
\usepackage{amsmath}
\usepackage{amsfonts}
\usepackage{amssymb}
\usepackage[final]{graphicx}
\usepackage{bm}
\usepackage[numbers]{natbib}

\graphicspath{{figures/}}

\author{T.~A.~Driscoll$^{1,3}$ \and R.~J.~Braun$^{1,4}$ \and J.~K.~Brosch$^{1,2}$}
\title{Simulation of Parabolic Flow on an Eye-Shaped Domain with Moving Boundary}
\institute{
$^1$Department of Mathematical Sciences University of Delaware, Newark, Delaware, 19716. \\
$^2$Current address: Arotech, 1229 Oak Valley Drive, Ann Arbor, Michigan 48108.\\
$^3$Email: driscoll@udel.edu; http://orcid.org/0000-0002-1490-2545.\\
$^4$http://orcid.org/0000-0002-5941-4166.
}

\usepackage{hyperref}
\hypersetup{
    colorlinks=true, 
    linktoc=all,     
    linkcolor=blue,  
}

\newcommand{\aref}[1]{\hyperref[#1]{Appendix~\ref{#1}}}

\def\equationautorefname~#1\null{(#1)\null}

\newcommand{\compdom}{\mathcal{C}}
\newcommand{\stripdom}{\mathcal{R}(t)}
\newcommand{\eyedom}{\mathcal{E}(t)}
\newcommand{\xhat}{\hat{x}}
\newcommand{\yhat}{\hat{y}}
\newcommand{\phihat}{\hat{h}}
\newcommand{\xtil}{\tilde{x}}
\newcommand{\ytil}{\tilde{y}}

\newcommand{\dd}[2]{\frac{d #1}{d #2}}
\newcommand{\pp}[2]{\frac{\partial #1}{\partial #2}}
\newcommand{\bfq}{\mathbf{q}}
\newcommand{\bfn}{\mathbf{n}}
\newcommand{\bfv}{\mathbf{v}}
\newcommand{\bfw}{\mathbf{w}}
\newcommand{\co}[1]{\operatorname{cx}\!\vphantom{\bigl(}\left(#1\right)}
\newcommand{\costar}[1]{\overline{\operatorname{cx}}\!\left(#1\right)}
\newcommand{\bigcdot}{\,\raisebox{-0.2ex}{\scalebox{1.4}{\ensuremath\cdot}}\,}
\newcommand{\nabtil}{\smash{\widetilde{\nabla}}}

\begin{document}

\maketitle

\begin{abstract}
    During the upstroke of a normal eye blink, the upper lid moves and paints a thin tear film over the exposed corneal and conjunctival surfaces. This thin tear film may be modeled by a nonlinear fourth-order PDE derived from lubrication theory. A challenge in the numerical simulation of this model is to include both the geometry of the eye and the movement of the eyelid. A pair of orthogonal and conformal maps transform a square into an approximate representation of the exposed ocular surface of a human eye.  A spectral collocation method on the square produces relatively efficient solutions on the eye-shaped domain via these maps. The method is demonstrated on linear and nonlinear second-order diffusion equations and shown to have excellent accuracy as measured pointwise or by conservation checks. Future work will use the method for thin-film equations on the same type of domain.
\end{abstract}

\begin{acknowledgements}
This work was supported by NSF grant DMS-1412085.  Any opinions, findings, and conclusions or recommendations
expressed in this material are those of the author(s) and do not necessarily
reflect the views of the National Science Foundation.
\end{acknowledgements}

\newpage
\section{Introduction}
\label{sec:intro}

The ocular tear film is critical for good vision and eye health.
In health, the tear film protects
the ocular surface with moisture, helps transport waste away from the ocular surface,
and provides a smooth optical surface for visual function \cite{DEWS}.
The tear film has multiple layers \cite{Mishima65,Ehl65,Norn79,BronTiffRev04,Gipson10}, 
but the thickest layer is composed primarily of water \cite{HollyLemp77}.
The pre-corneal tear film is the tear film located directly on the cornea, and that is what we refer to when we say tear film here.

The tear film is a total of a
few microns thick in the center of the cornea after a blink\cite{King-SmithFink04,KSetal06,WangEtal03} 
and has a considerably thicker meniscus (about $0.065$mm or more) around the lid margins
\cite{GoldingBruce97,Palakru07,JohnsonMurphy06,HarrisonBegley08}, where the tear film climbs the wettable part of the eyelids.
A relatively smooth and uniform film must be re-formed rapidly after each blink to enable vision with minimal interruption.   
Mathematical models have attempted to capture the dynamics of the tear film at various levels; recent reviews have appeared 
summarizing many efforts \cite{Braun12,BraunKing-Smith15}.
The term ``blink cycle'' is used to mean the combined periods of a single blink, in which the superior lid moves down toward the
inferior lid and then returns to its original position,
together with the interblink period separating two blinks.
Using one-dimensional models, a few papers have solved models for blink cycles \cite{BrKS07,Heryudono_Single_2007,ZubkovBreward12,Deng_A_2013,Deng_Heat_2014}, while others have solved for the combined
opening and interblink phases \cite{JonesPlease05,JonesMcElwain06,AydemirBreward10,MakiEtal08,JossicLefevre09,BruBrew14}.

The shape of the palpebral fissure (open eye shape) has also been incorporated into mathematical models of
tear film dynamics \cite{MakiBraun10a,MakiBraun10b}.
In these models, a piecewise polynomial boundary was created that was fit to a digital image of an open eye.
Besides specifying the tear film thickness at the boundary; they specified either the pressure
\cite{MakiBraun10a} or the flux \cite{MakiBraun10b} of the aqueous fluid.
Their simulations recovered features seen in previous 1D models and captured some experimental
observations of the tear film dynamics around the lid margins.  In \cite{MakiBraun10b},
the {\it in vivo} lacrimal supply and drainage mechanisms were simplified, and
they imposed a time-independent flux boundary condition; under some conditions,
they recovered flow around the outer canthus as seen experimentally.  The models were extended
to include important effects such as evaporation \cite{LiBraun14} and osmolarity transport
\cite{LiBraun16}.

To our knowledge, no models for a 2D blinking eye shape have been published.  There have been computational models for thin films on general surfaces that are time-independent \cite{GreBerSap06,RoyEtal02}, as well as analytical lubrication models on curved surfaces \cite{Howell03,BraunUsha12}.   Ultimately, we are interested in a blinking eye shape domain to study tear film dynamics, and so we focus here on the time dependent domain; there appears to be far less work done regarding moving domains.  In this work we introduce a model eye-shaped domain that blinks.  The upper and lower lids are arcs of circles and meet at a finite angle representing the canthi of the eye.  This domain approximates the palpebral fissure, but it has mathematical advantages over a closer approximation to the eye opening.  The eye-shaped domain can be conformally mapped to an infinite strip such that by making the upper side of the strip move vertically, the top edge of the eye shape can be made to move like the upper eyelid during a blink.  Futhermore, the infinite strip can be mapped to a square with edges at the fixed locations $x=\pm 1$ and $y=\pm 1$; this domain is convenient for applying Chebyshev spectral discretization of the spatial coordinates. In our models we solve the problems in this convenient computational domain and map the solution back to the eye-shaped domain for graphical representation of the results.

As a first step, we solve model problems that are second order in the spatial derivatives on this moving model eye domain.  We solve the (linear) heat equation, a version of a porous medium equation, and a nonlinear diffusion equation which mimics aspects of the thin fluid film equations. Nonlinear diffusion models in general, and porous medium equations in particular, have been of considerable interest (e.g., \cite{Kath82,Aronson1986,Witelski98}). We study only a small set of the possibilities here as proof of concept for use in nonlinear problems on this kind of moving domain.  The nonlinear models may have terms that have both positive and negative exponents, which imitates some aspects of the tear film models where the tear film is assumed to wet its substrate (the cornea in the case of \cite{WinAndBra10}).   We do not attempt to capture some aspects of possible solutions to the porous medium equations such as compact support \cite{Witelski98} or waiting times \cite{Kath82}.

We first state the test problems of interest in Section~\ref{sec:test-problems}.  We then set up the domain mappings to transform from the blinking eye shape to a fixed square computational domain (and back) in Section~\ref{sec:domains}.  The computational approach is discussed in Section~\ref{sec:comp-meth} and results are given in Section~\ref{sec:numresults}.  Finally, we discuss the results and future directions in Section~\ref{sec:Discussion}.

\section{Model problems}
\label{sec:test-problems}

We describe a solution method for the problem
\begin{equation}
  \label{eq:pde-div}
  h_t + \nabla \bigcdot \bfq(h,h_x,h_y) = 0, \qquad (x,y)\in \eyedom,
\end{equation}
where $h(t,x,y)$ is a dependent variable, meant to stand for the thickness of fluid, and $\bfq$ is a known flux function. The domain $\eyedom$, shown in \autoref{fig:mappings}, is an idealized eye shape whose top boundary moves in a prescribed fashion to be explained in section~\ref{sec:domains}.

\subsection{Linear diffusion}
The case $\bfq=-\kappa \nabla h$ for constant $\kappa$ results in the heat equation. It is of interest as the most basic second-order problem and because it has an explicit free-space solution that we can exploit to check the accuracy of our numerical method. Define the heat kernel
\begin{equation}
    \label{eq:heat-kernel}
    K(t,x,y) = \frac{1}{4\pi t} \exp{\left( -\frac{x^2+y^2}{4\kappa t} \right)}.
\end{equation}
Let $t_0>0$ and $(x_0,y_0)\in\mathcal{E}(0)$. If we use the initial condition $h(0,x,y)=K(t_0,x-x_0,y-y_0)$ and the Dirichlet condition
\begin{equation}
    \label{eq:heat-dirichlet}
    h(t,x,y) = K(t+t_0,x-x_0,y-y_0), \quad (x,y)\in \partial \eyedom,
\end{equation}
then equation~\autoref{eq:heat-kernel} holds throughout the domain for all $t$. This situation simulates the response to a point source at $(x_0,y_0)$ introduced at time $-t_0$.

\subsection{Nonlinear diffusion}

The nonlinear case $\bfq=-\psi(h) \nabla h$ can be made somewhat like the thin-film equation.  For example if $\psi(h)=\kappa(h^3+\beta h^{-3})$,
then the PDE resembles the case of a thin fluid film subject to van der Waals forces; this situation was studied in Winter \emph{et al} \cite{WinAndBra10} and many other places (e.g., \cite{JiWitelski17}).

Since exact solutions to this nonlinear problem are unknown to us, we use mass conservation as a check on the accuracy of our numerical solutions.  We create no-flux boundary conditions that keep the total mass of the solution fixed, and use deviation from that quantity as a proxy for the error in the
solution.  For any flux function $\mathbf{q}$, suppose that
\begin{equation}
    \label{eq:noflux-eye}
    \bfn\bigcdot \mathbf{q} - (\bfv\bigcdot \bfn) h = 0 \quad \text{on }\partial \eyedom,
\end{equation}
where $\bfn$ is the unit outward normal and $\bfv$ is the velocity of a point on the boundary. By the Reynolds Transport Theorem \cite{Acheson90} for the PDE~\eqref{eq:pde-div}, this condition ensures conservation of the total ``mass''
\begin{equation}
  \label{eq:mass-eye}
  M(t) = \int_{\eyedom} h(x,y)\, dA.
\end{equation}
If a boundary point is stationary, then~\eqref{eq:noflux-eye} is just a homogeneous Neumann condition. Equation~\eqref{eq:noflux-eye} can be generalized to prescribe any flux that varies in time and along the boundary.

\section{Domain and coordinate mappings}
\label{sec:domains}

\begin{figure}
        \begin{center}
            \includegraphics[width=\textwidth]{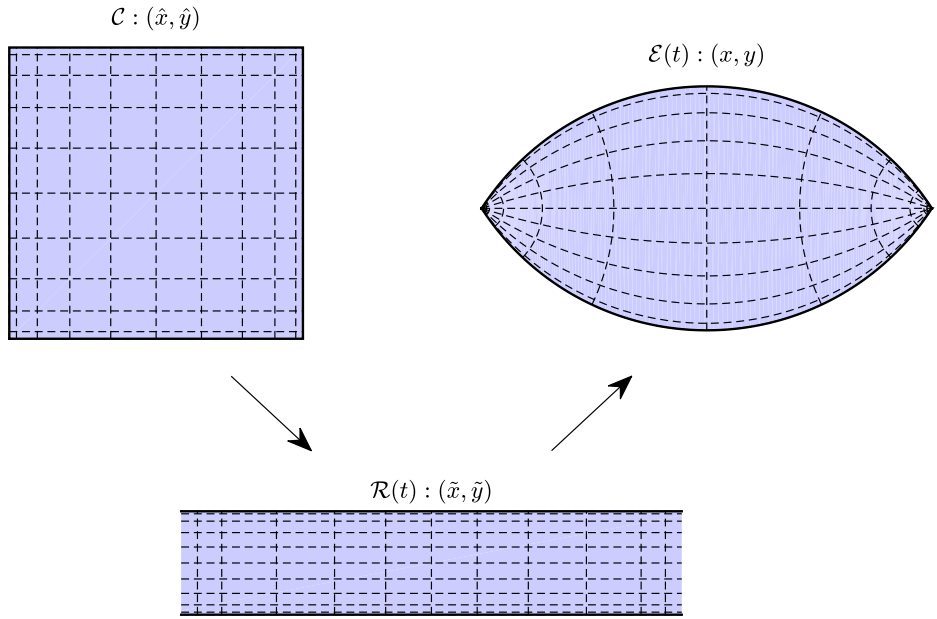}
        \end{center}
        \caption{Plot of the three domains used in this paper. The physical domain $\eyedom$ is meant to approximate the surface of the human eye, with an upper lid that moves in prescribed fashion. This region is mapped conformally using complex variables to an infinite strip $\stripdom$ whose top edge moves up and down as the eyelid moves in $\eyedom$. Finally, a time-varying map pulls the strip back to a fixed square $\compdom$, where a Chebyshev spectral method is applied for computations.}
    \label{fig:mappings}
\end{figure}

The PDE domain $\eyedom$ is determined by a pair of two-dimensional coordinate changes, as shown in \autoref{fig:mappings}.
Let $z=x+iy$ be the complex form of the coordinates in $\eyedom$, and let $\tilde{z}=\xtil+i\ytil$ be another complex coordinate. Define the infinte strip $\stripdom$ by
\begin{equation}
  \label{eq:stripdom}
  \stripdom = \{ \xtil+i\ytil\in\mathbb{C}: -1 < |\ytil| < \lambda(t) \},
\end{equation}
where $-1 < \lambda(t)\le 1$ is a prescribed function satisfying $\lambda(0)=1$ (which will correspond to a fully open eye). Then $\eyedom$ is the image of $\stripdom$ under the conformal map
\begin{equation}
  \label{eq:strip2eye}
    z = f(\tilde{z}) = \tanh\left( \frac{\tilde{z}}{2} \right),
\end{equation}
or equivalently,
\begin{equation}
  \label{eq:strip2eyexy}
  x(\xtil,\ytil) = \frac{\sinh(\xtil)}{\cos(\ytil) + \cosh(\xtil)} , \qquad
  y(\xtil,\ytil) = \frac{\sin(\ytil)}{\cos(\ytil) + \cosh(\xtil)}.
\end{equation}
The image of a line with fixed $\ytil$ is an arc ending at $(\pm1,0)$ of a circle centered at $(0,-\cot{\ytil})$~\cite[Appendix 2, Figure 20]{BrCh2004}. Two such arcs form the upper and lower boundaries of our eye-shaped region at all times. An alternative to our choice for a hyperbolic tangent map is to use parabolae to represent the edges of $\eyedom$; however, this choice leads to a much more complex and time dependent map \cite{Ivanov95}.  

Note that as $\xtil\to\pm\infty$, $(x,y)\to(\pm1,0)$, where $\eyedom$ has two corners whose angles change with time. As explained below, in practice we truncate $\stripdom$ by bounding $|\xtil|$ and therefore excise from $\eyedom$ exponentially small regions around the corners, leaving concave curves that always meet the upper and lower ``eyelids'' at right angles. The corners otherwise are likely to introduce singularities in the PDE that are of no interest in the tear film simulation problem and nontrivial to capture numerically.

A convenient feature of our setup is that the map $f$ between $\stripdom$ and $\eyedom$ is constant in time, yet the moving upper boundary of the strip $\stripdom$ maps to a curve that bears a plausible qualitative resemblance to a moving upper eyelid.  The geometry of $\stripdom$ and its moving boundary are easy to deal with computationally. Specifically, let us define an additional change of variables to a fixed computational domain $(\xhat,\yhat)\in \compdom=[-1,1]^2$:
\begin{equation}
    \label{eq:comp2strip}
    \xtil = \frac{\gamma\xhat}{\alpha^2 - \xhat^2}, \qquad \ytil = \frac{1}{2}(\yhat+1)(\lambda(t) + 1) - 1,
\end{equation}
for some $\gamma>0$ and $\alpha\ge 1$. If $\alpha=1$ then the image of $\compdom$ at time $t$ is all of $\stripdom$, but in practice we choose $\alpha > 1$ so that $|\xtil|\le \gamma(\alpha^2-1)$. The nonlinearity of the map between $\xhat$ and $\xtil$ allows us to compensate for the tendency of $f$ to crowd points near the corners of $\eyedom$.

\subsection{PDE transformation}

The change of variables from $\stripdom$ to $\eyedom$ causes transformation of the PDE \eqref{eq:pde-div} due to the chain rule. The most compact and elegant way to express these effects is by exploiting a correspondence between planar calculus and complex variables, as described in the appendix.

The additional change of variables from $\compdom$ to $\stripdom$ requires another application of the chain rule. The time dependence of the map between $y$-coordinates in~\eqref{eq:comp2strip} has a significant consequence. It must be the case that
\begin{equation}
  \label{eq:chainrule}
  \begin{split}
    \frac{\partial}{\partial{\ytil}} &=
    \frac{\partial \yhat}{\partial \ytil} \, \frac{\partial}{\partial{\yhat}}\\
    \frac{\partial}{\partial t} &=
    \frac{\partial \yhat}{\partial t} \, \frac{\partial }{\partial{\yhat}} + \frac{\partial}{\partial t}
  \end{split}.
\end{equation}
Hence, if we define $\hat{h}(\xhat,\yhat)=\tilde{h}(\xtil,\ytil)=h(x,y)$, then
\begin{equation}
  \label{eq:pde-chainrule}
   \hat{h}_t  = \tilde{h}_t - \pp{\yhat}{t}\,\hat{h}_{\yhat}= \tilde{h}_t  + \frac{\dot{\lambda}(1+\yhat)}{\lambda+1}\,\hat{h}_{\yhat},
\end{equation}
where $\dot{\lambda}=d\lambda/dt$. Thus an extra term is added to whatever expression is computed for $\tilde{h}_t$ in order to pass to $\hat{h}_t$, which is then evolved in the computational domain.

\subsection{No-flux condition}
\label{sec:no-flux-condition}

The no-flux boundary condition~\eqref{eq:noflux-eye} also transforms when mapping to $\stripdom$. We use the definitions and identities of the appendix to derive the result. Let $\bfq$, $\bfn$, and $\bfv$ denote the flux, normal, and velocity vectors in the original domain $\eyedom$. In complex terms,
\begin{equation*}
  \co{\bfq} = \psi(h) \co{\nabla h}
  = \psi(h) \frac{\co{\nabtil \tilde{h}}}{\overline{f'(\tilde{z})}}
  = \frac{\tilde{\bfq}}{\overline{f'(\tilde{z})}},
\end{equation*}
where $\tilde{\bfq}$ represents flux computed with respect to the variables in $\stripdom$. Since complex quantities transform by a factor of $f'$ under a conformal map, one finds
\begin{equation*}
  \co{\bfv} = f'(\tilde{z})\, \co{\tilde{\bfv}},
  \qquad
  \co{\bfn} = \frac{f'(\tilde{z})}{|f'(\tilde{z})|} \co{\tilde{\bfn}},
\end{equation*}
where $\tilde{\bfv}$ and $\tilde{\bfn}$ are velocity and normal vectors in the strip domain $\stripdom$. The no-flux condition~\eqref{eq:noflux-eye} becomes
\begin{align}
    0 = \bfn\bigcdot \mathbf{q} - (\bfn \bigcdot \bfv) h
      &= \operatorname{Re} \left[
        \co{\bfn}\costar{\bfq} - \co{\bfn}\costar{\bfv}h
        \right] \notag\\
      &= \operatorname{Re} \left\{
        \frac{f'(\tilde{z})}{|f'(\tilde{z})|}\co{\tilde{\bfn}}   \left[
        \frac{\costar{\tilde{\bfq}}}{{f'(\tilde{z})}}
        - \overline{f'(\tilde{z})}\, \costar{\tilde{\bfv}}h  \right]
        \right\} \notag\\
    &=  \frac{\tilde{\bfn}\bigcdot \tilde{\bfq}}{|f'(\tilde{z})|} - |f'(\tilde{z})| (\tilde{\bfn} \bigcdot \tilde{\bfv}) h.     \label{eq:noflux-strip}
\end{align}
The condition~\eqref{eq:noflux-strip} is much more straightforward to compute in $\stripdom$ than~\eqref{eq:noflux-eye} is in $\eyedom$, primarily because the normal and velocity vectors are very simple, e.g., $\tilde{\bfn}=\langle 0,1\rangle$ and $\tilde{\bfv}=\langle 0, \dot{\lambda}\rangle $ on the moving boundary.

\subsection{Mass}
\label{sec:mass}

The mass integral~\eqref{eq:mass-eye} transforms under coordinate changes as well. On the strip domain we have
\begin{equation}
    \label{eq:mass-strip}
    M(t) = \int_{\stripdom} h(\xtil,\ytil) J_\mathcal{R}(\xtil,\ytil) \,d\xtil\,d\ytil,
\end{equation}
where $J_\mathcal{R}$ is the Jacobian determinant.
Because the map $f$ from $\stripdom$ to $\eyedom$ is conformal, we have
\begin{equation}
  \label{eq:jacobian-strip}
  J_\mathcal{R}(\xtil,\ytil) = |f'(\tilde{z})|^{-2} = \bigl[\cosh(\xtil)+\cos(\ytil)\bigr]^{2}.
\end{equation}
Similarly, the integral may be computed in $\compdom$ by a product with another Jacobian:
\begin{equation}
\label{eq:mass-comp}
    M(t) = \int_{\compdom} h(\xhat,\yhat) J_\mathcal{R}(\xtil,\ytil) J_\compdom(\xhat,\yhat) \,d\xhat\,d\yhat,
\end{equation}
where
\begin{equation}
\label{eq:jacobian-comp}
    J_\compdom(\xhat,\yhat) = \frac{\gamma \xhat(1+\lambda(t))}{2\bigl(\alpha^2 - \xhat^2\bigr)^2}.
\end{equation}

\section{Computational method}
\label{sec:comp-meth}

We use a Chebyshev tensor-product spectral collocation discretization~\cite{SMiM}. Let $N$ be a discretization size and define the (one-dimensional) 2nd-kind Chebyshev points as
\begin{equation}
\label{eq:cheb-points}
x_k = -\cos \left( \frac{(k-1)\pi}{N-1} \right), \quad k=1,\ldots,N.
\end{equation}
A function $h(x)$ is represented by the $N$-vector $\bigl[h(x_1), \ldots, h(x_N) \bigr]^T$. These $N$ values define a unique polynomial interpolant $p(x)$ of degree less than $N$. Given analyticity of $h(x)$ on a complex region containing $[-1,1]$, the max-norm error $\|h-p\|_\infty$ converges as $O(C^{N})$ for some $C<1$, and  $\bigl[p^{(m)}(x_1), \ldots, p^{(m)}(x_N) \bigr]^T$ converge spectrally to sampled values of $h^{(m)}(x)$ for any integer $m>0$~\cite{ATAP}. The $N\times N$ matrix $D$ in the mapping
\begin{equation}
\label{eq:cheb-dm}
\begin{bmatrix}
p'(x_1) \\ \vdots \\  p'(x_N)
\end{bmatrix} =
D
\begin{bmatrix}
p(x_1) \\ \vdots \\  p(x_N)
\end{bmatrix}
\end{equation}
is the Chebyshev differentiation matrix. It is the discrete surrogate for the differentiation operator. We also use the analog of definite integration, known as Clenshaw--Curtis quadrature~\cite{SMiM}, in which one computes an $N$-vector $\bfw$ of quadrature weights such that
\begin{equation}
\label{eq:clencurt}
\bfw^T
\begin{bmatrix}
h(x_1) \\ h(x_2) \\ \vdots \\ h(x_N)
\end{bmatrix}
= \int_{-1}^{1}  p(x)\, dx \approx \int_{-1}^{1} f(x)\, dx.
\end{equation}

In the derivations that follow we use the $\odot$ operator, which is the Hadamard (elementwise) product between matrices (i.e., the \verb!.*! operator in MATLAB). We also find it convenient to make one more definition. Given any function $\phi(x)$, let
\newcommand{\diag}[1]{\operatorname{diag}_{#1}}
\begin{equation}
\label{eq:diag}
\diag{x}(\phi) =
\begin{bmatrix}
\phi(x_1) & & & \\
& \phi(x_2) & & \\
& & \ddots & \\
& & & \phi(x_N)
\end{bmatrix}.
\end{equation}
Left-multiplication by $\diag{x}(\phi)$ corresponds to pointwise multiplication by $\phi(x)$ on the grid.

\subsection{Method of lines}

Let $\xhat_i$, $i=1,\ldots,N_x$ and $\yhat_j$, $j=1,\ldots,N_y$ be two Chebyshev grids. To discretize the square domain $\compdom$, we use the tensor product grid
\begin{equation}
\label{eq:cheb-grid}
(\xhat_i,\yhat_j), \quad i=1,\ldots,N_x, \: j = 1,\ldots,N_y.
\end{equation}
We can represent the samples of a smooth function $\phihat$ on this grid as an $N_x\times N_y$ matrix satisfying $H_{ij}=\phihat(\xhat_i,\yhat_j)$. Suppose that $\widehat{D}_x$ and $\widehat{D}_y$ are appropriately sized differentiation matrices. Then the spatial derivatives of $\phihat$ on the grid are spectrally approximated by
\begin{equation}
\label{eq:cheb-par-diff}
\begin{split}
\left[ \pp{\phihat}{\xhat}(\xhat_i,\yhat_j) \right]_{i,j} &\approx \widehat{D}_x H \\
\left[ \pp{\phihat}{\yhat}(\xhat_i,\yhat_j) \right]_{i,j} \approx (\widehat{D}_y H^T)^T &= H \widehat{D}_y^T.
\end{split}
\end{equation}


The next step is to transform derivatives in $\compdom$ to those in the strip $\stripdom$. Considering~\eqref{eq:comp2strip} and~\eqref{eq:chainrule}, we define the corresponding one-dimensional operators
\begin{subequations}
  \label{eq:diffmat-strip}
  \begin{align}
    \widetilde{D}_x &=  \diag{\xhat} \left[ \frac{\left(\alpha ^2-x^2\right)^2}{\gamma (\alpha ^2+x^2)}\right]  \widehat{D}_x\\
    \widetilde{D}_y(t) &= \frac{2}{\lambda(t)+1} \widehat{D}_y.
  \end{align}
\end{subequations}
Then left-multiplication by $\widehat{D}_x$ and right-multiplication by $\widehat{D}_y^T$ perform partial derivatives on the grid in $\stripdom$, as in~\eqref{eq:cheb-par-diff}.


Next we consider the effect of the transformation between $\stripdom$ and the physical domain $\eyedom$. This is most conveniently done by exploiting the connection between planar vector calculus and complex variables, as described in the appendix. Define the matrices
\begin{align}
  \label{eq:heat-mol-fprime}
  F_{jk} &= \frac{1}{2}\operatorname{sech}^2[(\xtil_j+i\ytil_k)/2], \\
  E_{jk} &= 2\cosh^2[(\xtil_j+i\ytil_k)/2], \quad j=1,\ldots,N_x, \: k = 1,\ldots,N_y,
\end{align}
in which $(\xtil_j,\ytil_k)$ are the images of the grid points in $\compdom$ obtained from~\eqref{eq:comp2strip}. The entries of $F$ are values of $f'$ at the grid points, and the entries of $E$ are the values of the complex derivative of the inverse map $f^{-1}$. Also define a matrix $\Psi$ whose $(j,k)$ entry is $\psi(H_{jk})$. Then by~\eqref{eq:cxcalc1} and~\eqref{eq:cxxform1}, $\nabla h$ is computed from $\widetilde{\nabla}h$ through pointwise division by $f'$. This leads to
\begin{align}
  \label{eq:grid-grad-tilde}
  \widetilde{G} &= \widetilde{D}_x H + i H \widetilde{D}_y(t)^T, \\
  \label{eq:grid-grad}
  G &= \overline{E} \odot \widetilde{G},\\
  \label{eq:grid-flux}
  Q &= -\Psi \odot G,
\end{align}
as the grid approximations to $\nabtil h$, $\nabla h$, and $\bfq=-\psi(h)\nabla h$; note that the real and imaginary parts of $G$ correspond to $h_x$ and $h_y$ respectively, and the overline denotes complex conjugation. In similar fashion we use~\eqref{eq:cxcalc2} and~\eqref{eq:cxxform2} to derive the following discretization of $\nabla\bigcdot \bfq$ on the grid:
\begin{equation}
  \label{eq:grid-div}
  V = \operatorname{Re} \left\{
          E \odot \Bigl[ \widetilde{D}_x Q - i Q \widetilde{D}_y(t)^T \Bigr]
    \right\}.
\end{equation}

Finally, recall from~\eqref{eq:pde-chainrule} that a $t$-derivative in $\stripdom$ introduces an additional term in $\compdom$.
Hence the discretization of the PDE in $\compdom$ can be expressed as the matrix ODE
\begin{equation}
  \label{eq:grid-pde}
  \dd{H}{t}  = \Phi(H) :=  V  +
 \frac{\dot{\lambda}(t)}{\lambda(t)+1} \diag{\yhat}\bigl({1+\yhat}\bigr) H \widehat{D}_y^T.
\end{equation}
We should expect this system to be rather stiff. An $N\times N$ Chebyshev differentiation matrix has spectral radius $O(N^2)$, which becomes $O(N^4)$ when applied twice in each dimension.

\subsection{Boundary conditions}

To impose boundary conditions we modify the ODE system~\eqref{eq:grid-pde} by replacing differential equations at the boundary nodes by algebraic expressions of the appropriate conditions.
Let $B$ and $B'$ be the $N_x\times N_y$ \emph{boundary indicator} matrices
\begin{equation}
  \label{eq:boundary-mask}
  B = \begin{bmatrix}
    1 & 1 & \cdots & 1 & 1 \\
    1 & 0 & \cdots & 0 & 1 \\
    \vdots & \vdots &  & \vdots & \vdots \\
    1 & 0 & \cdots & 0 & 1 \\
    1 & 1 & \cdots & 1 & 1
  \end{bmatrix}, \qquad
  B' = \begin{bmatrix}
    0 & 0 & \cdots & 0 & 0 \\
    0 & 1 & \cdots & 1 & 0 \\
    \vdots & \vdots &  & \vdots & \vdots \\
    0 & 1 & \cdots & 1 & 0 \\
    0 & 0 & \cdots & 0 & 0
  \end{bmatrix}.
\end{equation}
For any boundary condition, we replace~\eqref{eq:grid-pde} by
\begin{equation}
\label{eq:boundary-dirichlet}
{B}' \odot H_t = B \odot R + {B}' \odot \Phi(H),
\end{equation}
where $R$ is the residual value of boundary condition equations at the boundary points. For instance, if $W(r,s,t)$ is an $N_x\times N_y$ matrix whose boundary entries represent Dirichlet values for the solution, then $R=H-W$. Equation~\eqref{eq:boundary-dirichlet} is an index-1 differential--algebraic equation (DAE) for $H$. We use the MATLAB function \texttt{ode15s} to solve such systems.

For no-flux conditions, the boundary entries of $R$ are computed from the condition~\eqref{eq:noflux-strip}, which is equivalent to~\eqref{eq:noflux-eye}. For the left and right boundaries of $\stripdom$, i.e. the first and last rows on the grid, the values are $\operatorname{Re}(\widetilde{G})$, and for the bottom of $\stripdom$ (first column of the grid), the values are $\operatorname{Im}(\widetilde{G})$. Finally, for the moving top boundary, the values come from
\begin{equation*}
  E \odot \Psi \odot \operatorname{Im}(\widetilde{G}) - \dot{\lambda}(t) F\odot H.
\end{equation*}

\subsection{Computation of mass}
\label{sec:computation-of-mass}

The mass integral $M(t)$ in~\eqref{eq:mass-eye} is computed using Clenshaw--Curtis quadrature as defined in~\eqref{eq:clencurt}.  We will need to incorporate the Jacobians from the mappings as indicated in~\eqref{eq:mass-comp}. Let $\hat{\bfw}_x$ and $\hat{\bfw}_y$ be the Clenshaw--Curtis weight vectors of length $N_x$ and $N_y$ respectively. We represent the Jacobians in the integrand in~\eqref{eq:mass-comp} by the matrix $S$, where
\begin{equation*}
S_{ij}(t) = J_\mathcal{R}(\xtil_i,\ytil_j) J_\compdom(\xhat_i,\yhat_j), \qquad i=1,\ldots,N_x,\quad j=1,\ldots,N_y.
\end{equation*}
Then the mass integral is computed as
\begin{equation}
\label{eq:mass-clencurt}
M(t) \approx \hat{\bfw}_x^T \bigl[S(t) \odot H(t)\bigr] \hat{\bfw}_y.
\end{equation}

\subsection{Lid motion}

The observed lid motion during a blink is primarily by the upper lid descending to the lower lid, and the outer canthus being pulled in
about 10-20\% of the horizontal width of the palpebral fissure \cite{Doane_Interaction_1980}.  We neglect the horizontal contraction here, and move only the upper boundary while keeping the locations corresponding to the canthi fixed.   Various formulations for the motion of the upper lid have been proposed.  Berke and M\"uller \cite{BeMu98} proposed a product of a monomial with decaying exponential for the opening phase of the blink; a similar function was proposed in \cite{JossicLefevre09}.   A more complex function was proposed in \cite{AydemirBreward10}, which has been used in several subsequent papers.  For full blink cycles, a sinusoidal function was used in \cite{BrKS07} and a more realistic function was proposed in \cite{Heryudono_Single_2007}, which has also been used subsequently.  Here we use a simplified periodic lid motion function which improves upon the sinusoidal version, but is still simpler than the realistic versions.

\begin{figure}
  \centering
  \includegraphics{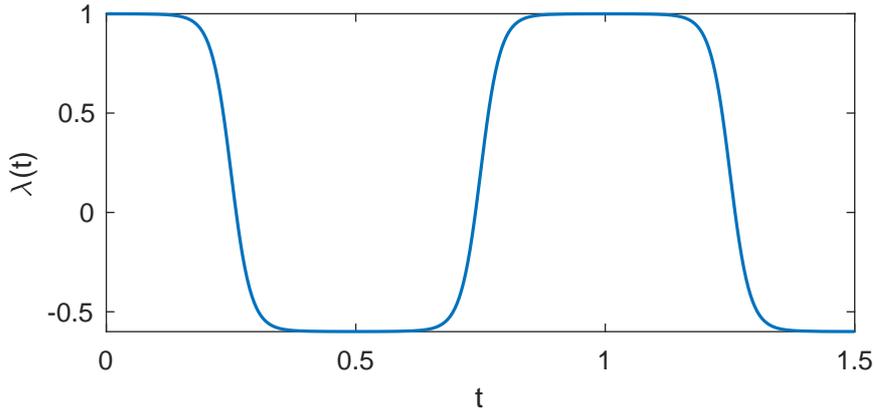}
  \caption{One and a half periods of the upper lid motion function~\autoref{eq:rho}, with $c=0.8$ and $\nu=1$.}
  \label{fig:lidmotion}
\end{figure}

In all of the experiments reported below, we use the lid motion function
\begin{equation}
   \label{eq:rho}
   \lambda(t) = 1-c + c \tanh \bigl( 4 \cos(2\pi \nu t) \bigr),
\end{equation}
where $c$ represents the maximum percentage closure of the eye and $\nu$ is the frequency of the blink. This function is periodic with significant pauses at the fully open and fully closed phases, as is shown in \autoref{fig:lidmotion} for $c=0.8$ and $\nu=1$. This motion is prescribed in the strip domain $\stripdom$ and mapped to the lens domain via~\autoref{eq:strip2eyexy}.

\section{Numerical results}
\label{sec:numresults}

Computations were executed in MATLAB 2016b using \texttt{ode15s} as the time integrator for the DAE~\autoref{eq:boundary-dirichlet}. The reported computation times are for a 2014 iMac (4 GHz Intel Core i7) running macOS 10.12.2.

\subsection{Heat equation with known solution}
\label{sec:heat-known}

The first results are for~\eqref{eq:pde-div} with $\bfq=-\nabla h$, i.e. the heat equation with unit diffusion constant. For the lid motion we used $c=0.8$ as the closure fraction and $\nu=16$ as the frequency in~\eqref{eq:rho}. The initial and boundary conditions were used to make the exact solution equal the free-space heat kernel as in~\eqref{eq:heat-dirichlet} with $t_0=0.01$, $x_0=0.1$, and $y_0=0.2$.
Using a 28$\times$24 grid and time integration error tolerance set to $10^{-9}$, the solution took 11.3 seconds to compute two full blink cycles (up to $t=1/8$).

\begin{figure}
  \centering
  \includegraphics{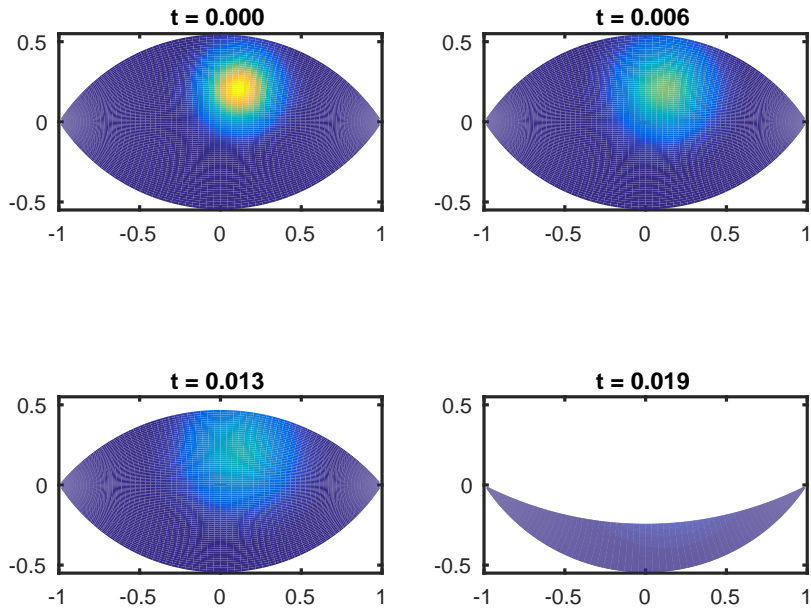}
  \caption{Computed solution of the heat equation with a free-space Dirichlet boundary condition. The eye is initially fully open with a Gaussian solution ranging from zero (dark blue) to 8 (yellow). The Gaussian spreads as the eyelid moves, and it is unaffected by the moving boundary.}
  \label{fig:heatfreeblink}
\end{figure}

\begin{figure}
  \centering
  \includegraphics{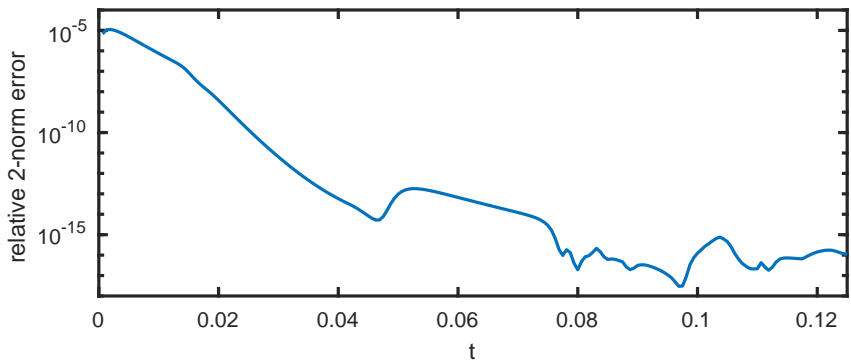}
  \caption{Error in the computed solution of \autoref{fig:heatfreeblink}, as measured by square integration over the physical domain and relative to the exact solution.}
  \label{fig:heatfreeerror}
\end{figure}

Snapshots of the computed solution are shown in \autoref{fig:heatfreeblink}. The computed solution shows two expected tendencies: there are no effects from the moving boundaries, and the solution becomes smoother with time in a symmetric fashion around $(x_0,y_0)$.  The decay of the free spaced Green's function is relatively fast compared to the lid motion in this case.  In \autoref{fig:heatfreeerror} we plot the relative error
\begin{equation*}
  \frac{\| h_{\text{comp}}(x,y,t)-h(x,y,t)\|_{\eyedom}}{\|h(x,y,t)\|_{\eyedom}} ,
\end{equation*}
in which $h$ is the known exact solution and the norm is defined using square integration over the physical domain (as performed spectrally in $\compdom$). We see that the solution is accurate to five digits initially and improves at more or less an exponential rate as the solution proceeds, until the solution is zero relative to the initial condition.  This level of accuracy is expected for such a smooth solution using spectral methods.

\subsection{Porous medium equation}
\label{sec:porous}

Our next example is the porous medium equation with flux
\begin{equation}
  \label{eq:porous}
  \bfq=-[(1-\kappa)h + \kappa]\nabla h,
\end{equation}
for $0<\kappa\le 1$. The lid motion function was \autoref{eq:rho} with closure fraction $c=0.7$ and frequency $\nu=1$.  The initial condition was specified in the strip domain $\mathcal{R}(0)$ as
\begin{equation}
  \label{eq:initcond}
   \tilde{h}(0,\tilde{x},\tilde{y}) =  1 - 0.8 e^{-6(\tilde{y}+0.2)^2-4(\tilde{x}-1)^2}.
\end{equation}
This is not a Gaussian function in the physical domain $\mathcal{E}(0)$. We chose the mass-conserving flux boundary condition as described in \autoref{sec:no-flux-condition}.

\begin{figure}
  \centering
  \includegraphics{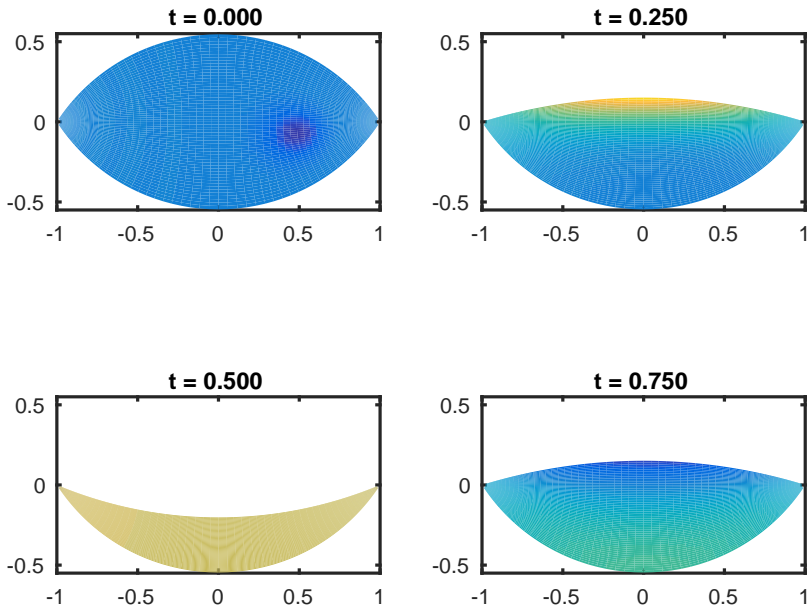}
  \caption{Computed solution of the porous medium equation~\autoref{eq:porous} with a mass-conserving boundary condition. The solution ranges from about zero (dark blue) to about 4 (near $t=0.25$, yellow).}
  \label{fig:porousblink}
\end{figure}

For this problem we used a $32\times 48$ grid and a time integrator tolerance of $10^{-9}$. It took approximately 73 seconds to solve for two complete blink cycles with $\kappa=0.5$; the solution at selected times is presented in \autoref{fig:porousblink}. (Computation times for other values of $\kappa$ were similar.)  For our parameters, the influence of the initial condition disappears from the solution quickly. The solution develops a boundary layer along the descending upper boundary where $h$ is increased; see $t=0.25$.  When the upper boundary is ascending, there is also a boundary layer but at these times $h$ is decreased near the boundary ($t=0.75$).  There is a pause in the lid position around $t=0.5$, which allows $h$ to become uniform across the domain.

\begin{figure}
  \centering
  \includegraphics{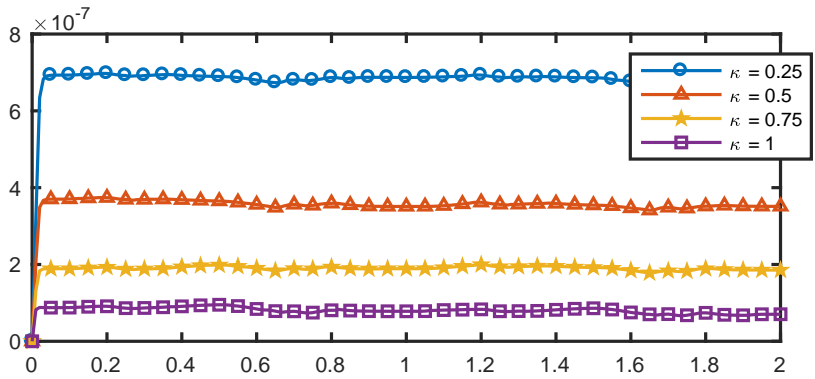}
  \caption{Mass change in computed solutions of the porous medium equation~\autoref{eq:porous}. The exact solution has constant total mass, so the size of this quantity is indicative of pointwise error in the numerical solution.}
  \label{fig:porousmass}
\end{figure}

Since the boundary conditions conserve mass while the numerical methods do not explicitly do so, we use the computed relative change in the total mass as an indicator of the accuracy of the computation. As demonstrated by the conservation results in \autoref{fig:porousmass} for several values of $\kappa$, we feel justified in claiming six-digit pointwise accuracy for the computation.  The size of the error appears to be related to the relative strength of the nonlinearity and changes very little after just a few time steps.  The relatively modest number of grid points appears to resolve the boundary layers and boundary conditions quite well.

\subsection{Thin-film analog}
\label{sec:film}

Our closest second-order analog to the thin-film problem is the equation
\begin{equation}
  \label{eq:film}
  \bfq=-[A-Bh^{-3}]\nabla h,
\end{equation}
for constants $A$ and $B$. The nonlinear term creates a preference for a flat solution at the value $(B/A)^{1/3}$. Once such a solution is established in one region, it can serve as a barrier to diffusion through that region. Our simulations used $A=1$ and $B=10^{-9}$. The lid motion had closure fraction $c=0.8$ and frequency $\nu=1$. The initial condition was constant: $h(0,x,y)=0.1$. The boundary condition again used the flux to conserve total mass in the exact solution.

\begin{figure}
  \centering
  \includegraphics{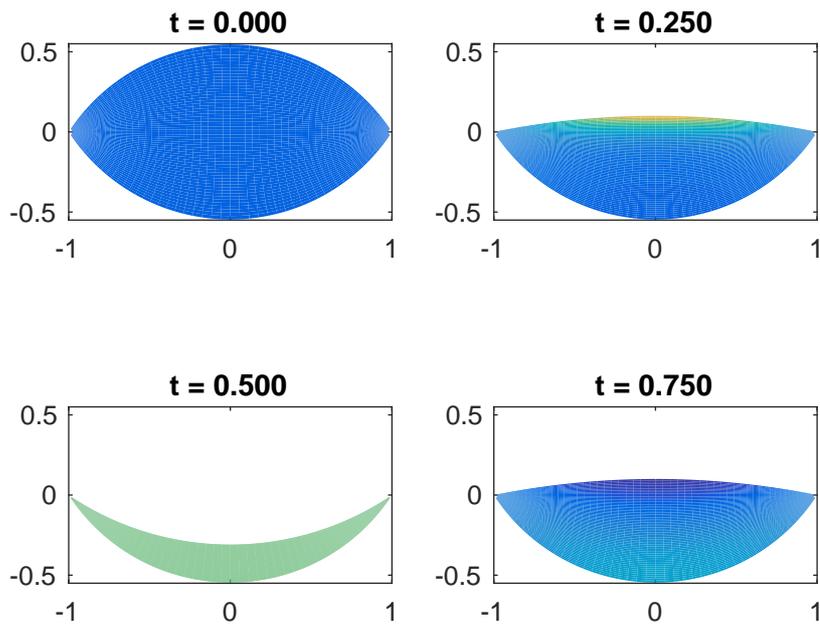}
  \caption{Computed solution of the thin-film analog equation~\autoref{eq:film} with a mass-conserving boundary condition. The solution ranges from zero (dark blue) to about 1.87 (yellow).}
  \label{fig:filmblink}
\end{figure}

\begin{figure}
  \centering
  \includegraphics{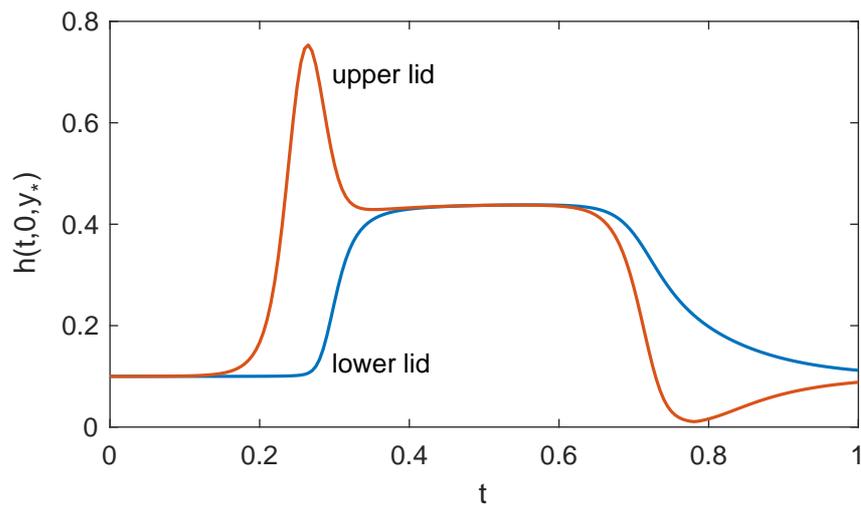}
  \caption{Film height at the centers of the upper and lower lids for computed solutions of the thin-film analog equation~\autoref{eq:film}. See also \autoref{fig:filmblink}.}
  \label{fig:filmlid}
\end{figure}

\begin{figure}
  \centering
  \includegraphics{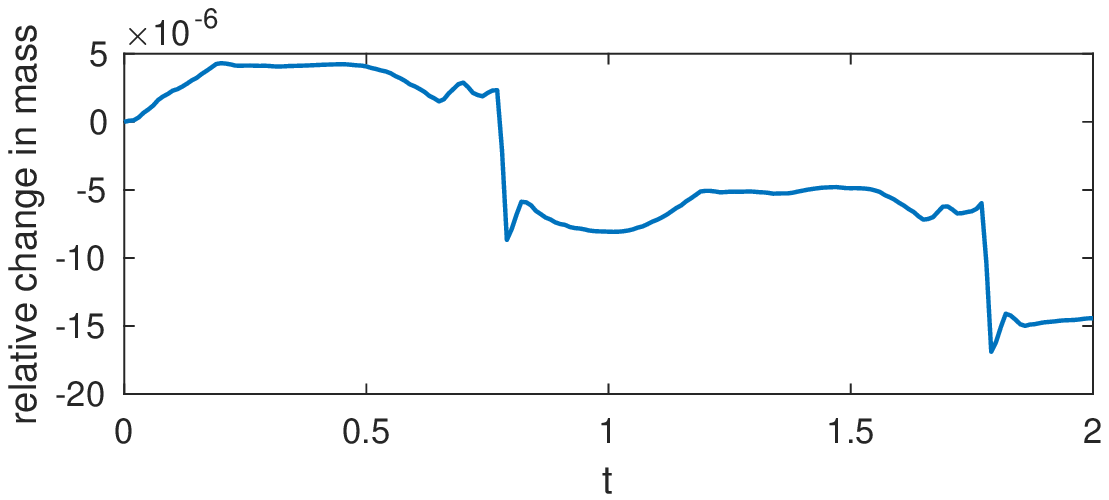}
  \caption{Mass change in computed solutions of the thin-film analog equation~\autoref{eq:film}. The exact solution has zero mass change.}
  \label{fig:filmmass}
\end{figure}

Using a $31\times 40$ grid with time integrator tolerance set to $10^{-7}$, computing up to $t=2$ (two blink cycles) took 40 seconds. Snapshots of the result are shown in \autoref{fig:filmblink}.  During the downward motion of the upper edge at $t=0.25$, the solution $h$ builds up at the boundary forming a steep local boundary layer around $y=0$.  The pause at $t=0.5$, at which the domain has minimum area, allows $h$ to become essentially uniform at about $h=0.42$, more than four times the value of the initial condition.  During the upward motion at $t=0.75$, $h$ is locally depleted near the center of the moving edge, again forming a boundary layer there.
\autoref{fig:filmlid} shows the solution at the midpoints (along $x=0$) of the upper and lower edges or ``lids'' of the moving domain.  The localized changes around the center of the moving upper boundary are clearly seen. Among the times $t_j=j/100$ for integer $j$, the solution at the upper lid takes a minimum value of $h\approx 0.01067$ at $t=0.78$ in the first cycle.  The tendency to quickly form a uniform state around $t=0.5$ is also clear.
The mass conservation error plotted in \autoref{fig:filmmass} suggests that we have obtained five-digit accuracy for the solution.

If we decrease the initial volume of the solution, the time integrator stagnates near the time of maximum upward velocity of the upper lid (corresponding to $t=0.78$ in the last case). This stagnation time seems to remain constant as the spatial grid is refined. We speculate that as the lid opens there may be insufficient material to prevent the solution from reaching zero, which appears to cause nonexistence of the PDE solution \cite{Aronson1986}.

\section{Discussion}
\label{sec:Discussion}

Conformal mappings allowed us to compute answers accurately in a simple computational square and transfer the results back to a moving eye-shaped domain.  
In particular, the map from the intermediate strip domain to the eye-shaped domain was time independent and its simplicity greatly facilitated the computation.  Also, the strip domain allowed a simple expression for the moving boundary.

Using a linear heat equation with an exact free-space Green's function solution, we demonstrated that we can accurately recover the solution to that problem.  We achieved high accuracy with a relatively small number of grid points in each direction, as is to be expected for a spectral method applied to a smooth solution.

We then moved to a type of porous medium equation in which the ``diffusivity'' was a linear function that interpolated between a porous medium equation and the linear heat equation.  We could again obtain solutions that conserved mass very well with relatively few grid points.  The conservation of mass was used a proxy for accuracy in the nonlinear equations in the absence of an exact solution.  For our final model problem, we studied a nonlinear diffusion equation with both positive and negative powers of the dependent variable, in an effort to mimic the expected dynamics of the tear film on the moving domain we used.  We observed transient boundary layers as may be expected for the thin film problem
(e.g., \cite{Heryudono_Single_2007,AydemirBreward10,Deng_A_2013,Deng_Heat_2014}).  The error in mass conservation was again small, with a relatively small number of grid points in each direction.  The numerical method did have difficulty when the thickness reached an expected equilibrium value; resolving that difficulty is beyond the scope of this paper.

Turning to future work, extending the model to treat the fourth order thin film problem is needed to more directly study the tear film.  Work on this extension is underway.

An important aspect of in vivo tear dynamics is the influx and efflux of tears during the blink cycle.
In his theory of the lacrimal system,
Doane\cite{Doane81} proposed that significant drainage along the lid margins {begins}
with the lids about halfway open, and that it {ends} up to 3s after the lids have fully opened.
New tear fluid is supplied from the lacrimal gland, which secretes the aqueous part of the tear film,
and
enters the exposed tear film from beneath the upper lid near the outer (temporal)
canthus \cite{Maurice73,HarrisonBegley08}.  The aqueous part of tears exits via the puncta,
which are small holes found near the nasal canthus.  Tear film models have incorporated lacrimal
gland supply and punctal drainage as well, in both
one-dimensional \cite{Heryudono_Single_2007,MakiEtal08,Deng_A_2013,Deng_Heat_2014}
and two-dimensional \cite{MakiBraun10b,LiBraun14} models.  In future work, we will incorporate the flux conditions of lacrimal gland and punctal
drainage as the previous 2D models, and compare the results with relevant experiments \cite{HarrisonBegley08,LiBraun14}.

Incorporating more realistic lid motion functions for the blink cycle \cite{Heryudono_Single_2007,ZubkovBreward12} would be a valuable extension as well.

\bibliographystyle{spbasic}
\bibliography{tearfilms-rjb}

\appendix
\section{Complex-variable expressions for 2D vector fields}

There are ways to express two-dimensional vector fields and vector calculus that can simplify derivations and computations relating to changing variables. For a vector field $\mathbf{r}(x,y)=\langle r_1(x,y),r_2(x,y) \rangle$, define
\begin{equation}
\label{eq:complexify}
\co{\mathbf{r}} = \co{\bigl\langle r_1(x,y),r_2(x,y)\bigr\rangle } = r_1(x,y) + ir_2(x,y)
\end{equation}
Then it is easy to confirm that for two vector fields $\mathbf{r}$ and $\mathbf{s}$,
\begin{equation}
\mathbf{r} \bigcdot \mathbf{s} = \operatorname{Re}\Bigl[\costar{\mathbf{r}}\co{\mathbf{s}}\Bigr],
\end{equation}
where the overline indicates complex conjugation of the result. Next define a gradient operator and its complexification by
\begin{equation}
  \label{eq:nabla}
  \nabla = \left\langle \frac{\partial}{\partial {x}}, \frac{\partial}{\partial {y}}\right\rangle \,,
      \quad \co{{\nabla}} = \frac{\partial}{\partial {x}} + i \frac{\partial}{\partial {y}},
      \quad \costar{{\nabla}} = \frac{\partial}{\partial {x}} - i \frac{\partial}{\partial {y}}\,.
\end{equation}
Then it is easily checked that for any differentiable scalar function $\phi$ and vector field $\mathbf{r}$,
\begin{subequations}
  \label{eq:cxvcalc}
  \begin{align}
    \co{\nabla \phi} &= \co{\nabla}\phi \label{eq:cxcalc1} \\
    \nabla \bigcdot \mathbf{r} &= \operatorname{Re}\Bigl[  \costar{\nabla} \co{\mathbf{r}} \Bigr] \label{eq:cxcalc2} \\
    \nabla^2 \phi &= \costar{\nabla}\co{\nabla} \phi. \label{eq:cxcalc3}
  \end{align}
\end{subequations}

Now suppose that we use a conformal map $f$ to change variables; specifically, let the connection between two planes be
\begin{equation}
  \label{eq:cxmap}
  z = x+iy = f(\xtil+i\ytil) = f(\tilde{z}).
\end{equation}
Let $\tilde{\phi}(\xtil,\ytil)=\phi(x,y)$ and $\tilde{\mathbf{r}}(\xtil,\ytil)=\mathbf{r}(x,y)$, subject to~\eqref{eq:cxmap}, and let~\eqref{eq:nabla} be extended to the case with tildes on $\nabla$, $x$, and $y$. Then we have the identities
\begin{subequations}
  \label{eq:cxxform}
  \begin{align}
    \co{\nabtil \tilde{\phi}} &= \overline{f'(\tilde{z})}\, \co{\nabla \phi}  \label{eq:cxxform1} \\
    \nabtil \bigcdot \tilde{\mathbf{r}} &= \operatorname{Re}\Bigl[  f'(\tilde{z}) \costar{\nabla} \co{\mathbf{r}} \Bigr] \label{eq:cxxform2} \\
    \nabtil^2 \tilde{\phi} &= \bigl|f'(\tilde{z})\bigr|^2\, \nabla^2 \phi,\label{eq:cxxform3}
  \end{align}
\end{subequations}
where the prime on $f$ indicates usual differentiation in a complex variable.

\end{document}